\documentclass[12pt]{article}

\usepackage{amsmath,amsthm,amsfonts,amssymb, color}

\newtheorem{theorem}{Theorem}[section]

\newtheorem{lemma}[theorem]{Lemma}

\newtheorem{remark}[theorem]{Remark}

\def\cA{\mathcal{A}}

\def\cE{\mathcal{E}}
\def\cF{\mathcal{F}}

\def\cH{\mathcal{H}}
\def\cI{\mathcal{I}}
\def\cS{\mathcal{S}}

\def\bR{\mathbb{R}}
\def\bC{\mathbb{C}}

\begin{document}

\title{Second order Lyapunov exponents for \\
parabolic and hyperbolic Anderson models}

\author{Raluca M. Balan\footnote{Department of Mathematics and Statistics, University of Ottawa,
585 King Edward Avenue, Ottawa, ON, K1N 6N5, Canada. E-mail
address: rbalan@uottawa.ca} \footnote{Research supported by a
grant from the Natural Sciences and Engineering Research Council
of Canada.}\and
Jian Song\footnote{Corresponding author. Department of Mathematics. University of Hong Kong, Hong Kong. E-mail address: txjsong@hku.hk}}

\date{April 5, 2017}
\maketitle

\begin{abstract}
\noindent In this article, we consider the hyperbolic and parabolic Anderson models in arbitrary space dimension $d$, with constant initial condition, driven by a Gaussian noise which is white in time. We consider two spatial covariance structures: (i) the Fourier transform of the spectral measure of the noise is a non-negative locally-integrable function; (ii) $d=1$ and the noise is a fractional Brownian motion in space with index $1/4<H<1/2$. In both cases, we show that there is striking similarity between the Laplace transforms of the second moment of the solutions to these two models. Building on this connection and the recent powerful results of \cite{huang-le-nualart15} for the parabolic model, we compute the second order (upper) Lyapunov exponent for the hyperbolic model. In case (i), when the spatial covariance of the noise is given by the Riesz kernel, we present a unified method for calculating the second order Lyapunov exponents for the two models.
\end{abstract}


\vspace{1mm}


\section{Introduction}

The goal of this article is to study the second order Lyapunov exponent
of the solution to the hyperbolic Anderson model:
\begin{equation}
\left\{\begin{array}{rcl}
\displaystyle \frac{\partial^2 u}{\partial t^2}(t,x) & = & \displaystyle \Delta u(t,x)+ u(t,x)\dot{W}(t,x), \quad t>0, x \in \bR^d \\[2ex]
\displaystyle u(0,x) & = &  1 , \quad x \in \bR^d\\[1ex]
\displaystyle \frac{\partial u}{\partial t}(0,x) & = & 0, \quad x \in \bR^d
\end{array}\right. \label{HAM} 
\end{equation}
driven by a zero-mean Gaussian noise $\dot{W}$ which is white in time and is spatially homogeneous with spatial covariance given by a tempered measure $\mu$ on $\bR^d$.
We consider also the parabolic Anderson model:
\begin{equation}
\left\{\begin{array}{rcl}
\displaystyle \frac{\partial u}{\partial t}(t,x) & = & \displaystyle \frac{1}{2}\Delta u(t,x)+ u(t,x)\dot{W}(t,x), \quad t>0, x \in \bR^d \\[2ex]
\displaystyle u(0,x) & = &  1 , \quad x \in \bR^d\\[1ex]
\end{array}\right. \label{PAM} 
\end{equation}
with the same noise $\dot{W}$ as above. We denote by $u^w$ and $u^h$ the solutions to equation \eqref{HAM}, respectively \eqref{PAM}. To simplify the writing, we use the convention that whenever a result holds for both equations \eqref{HAM} and \eqref{PAM}, we omit writing the indices $w$ and $h$.

The noise is given by an isonormal Gaussian process $W=\{W(\varphi); \varphi \in \cH\}$ with covariance $$E[W(\varphi)W(\psi)]=\langle \varphi,\psi \rangle_{\cH},$$ where $\cH$ is the completion of the space $C_0^{\infty}(\bR_{+} \times \bR^d)$ of infinitely differentiable functions with compact support on $\bR_{+} \times \bR^d$, with respect to the inner product:
\begin{equation}
\label{inner-product-H}
\langle \varphi,\psi \rangle_{\cH}
=\int_{0}^{\infty} \int_{\bR^d} \cF \varphi(t,\cdot)(\xi) \overline{\cF \psi(t,\cdot)(\xi)} \mu(d\xi) dt,
\end{equation}
and $\cF \varphi(t,\cdot)$ denotes the Fourier transform of the function $\varphi(t,\cdot)$. We define the Fourier transform of a function $\varphi \in L^1(\bR^d)$
by $\cF \varphi(\xi)=\int_{\bR^d}e^{-i \xi \cdot x}\varphi(x) dx$ for $\xi \in \bR^d$.
Here $\xi \cdot x$ is the inner product in $\bR^d$ and $|\cdot|$ is the Euclidean norm in $\bR^d$.

We are interested in two cases:
\begin{description}
\item[(i)] the Fourier transform of $\mu$ (in the space $\cS'(\bR^{d})$ of tempered distributions on $\bR^d$) is a non-negative locally integrable function $f$, i.e.
\begin{equation}
\label{Fourier-f-mu}
\int_{\bR^d}f(x)\varphi(x)dx=\int_{\bR^d}\cF \varphi(\xi)\mu(d\xi), \quad \mbox{for any} \quad \varphi \in \cS(\bR^d),
\end{equation}
     and the measure $\mu$ satisfies {\em Dalang's condition}:
    \begin{equation}
    \label{Dalang-cond}
    \int_{\bR^d}\frac{1}{1+|\xi|^2}\mu(d\xi)<\infty.
    \end{equation}
\item[(ii)] $d=1$ and $\mu(d\xi)=c_{H} |\xi|^{1-2H}d\xi$ with $1/4<H<1/2$ and $c_{H}=\Gamma(2H+1)\sin(\pi H)/(2\pi)$. In this case, the Fourier transform of $\mu$ in $\cS'(\bR)$ is the distribution $\Gamma$ defined by
    regularization: $\Gamma(\varphi)=H(2H-1)\int_{\bR}(\varphi(x)-\varphi(0))dx$  {(see e.g. p. 23-24 of \cite{gelfand-shilov64})}. We have:
    \begin{equation}
    \label{def-Gamma}
    \Gamma(\varphi)=\int_{\bR}\cF \varphi(\xi)\mu(d\xi), \quad \mbox{for any} \ \varphi \in \cS(\bR).
    \end{equation}
\end{description}

Existence (and uniqueness) of the solution to equation \eqref{HAM} has been proved recently in \cite{balan-song16} (for case (i), with a noise more general than here) and in \cite{BJQ} (for case (ii)).
In both cases (i) and (ii), it is also known that equation \eqref{PAM} has a unique solution (see \cite{huang-le-nualart15}).

Note that in case (i), the inner product in $\cH$ can be written as:
$$\langle \varphi, \psi \rangle_{\cH}=\int_{0}^{\infty}\int_{\bR^d} \int_{\bR^d} \varphi(t,x)\psi(t,y)f(x-y)dx dy dt.$$
An example for this case is the Riesz kernel
$f(x)=|x|^{-\alpha}$ with $0<\alpha<d$. For this example, relation \eqref{Fourier-f-mu} holds for $\mu(d\xi)=C_{d,\alpha}|\xi|^{-(d-\alpha)}d\xi$ with
\begin{equation}
\label{def-C-d-alpha}
 {
C_{d,\alpha}=\pi^{-d/2}2^{-\alpha}\,
\frac{\Gamma(\frac{d-\alpha}{2})}{\Gamma(\frac{\alpha}{2})},
}
\end{equation}
(see e.g. Lemma 1, p.117 of \cite{stein70}) and condition \eqref{Dalang-cond} holds if and only if
$\alpha<2$.

Whereas the moments of the solution to equation \eqref{PAM} have been studied intensively in the recent years (see e.g. \cite{dawson-salehi80}, \cite{FK13}, \cite{chen-hu-song-xing15} \cite{huang-le-nualart15}, \cite{chen17}), our results seem to be the first which give an exact calculation for the (upper) Lyapunov exponent of the solution to equation \eqref{HAM} in arbitrary space dimension $d$. The $p$-th order upper and lower Lyapunov exponents for the solution to equation \eqref{HAM} in dimension $d=3$ were studied in \cite{dalang-mueller09}, for case (i). The existence of the solution to the stochastic wave equation in arbitrary dimension $d$ was established in \cite{conus-dalang08}, again in case (i).


By definition, a {\it mild} solution of equation \eqref{HAM}, respectively \eqref{PAM}, is a square-integrable process $u$, which is adapted to the filtration induced by the noise $W$ and satisfies the integral equation:
\begin{equation}
\label{integral-eq}
u(t,x)=1+\int_{0}^{t}\int_{\bR^d}G(t-s,x-y)u(s,y)W(ds,dy),
\end{equation}
for any $t>0$ and $x\in \bR^d$, where $G=G^w$ is the fundamental solution of the wave equation, respectively $G=G^h$ is the fundamental solution of the heat equation. As in \cite{balan-song16}, the stochastic integral in \eqref{integral-eq} is interpreted in the Skorohod sense, but it can be shown to coincide with the It\^o-type integral considered in \cite{dalang99}.

We note in passing that $G^h(t,\cdot)$ is a rapidly decreasing function in $\bR^d$, whereas $G^w(t,\cdot)$ is an integrable function for $d=1$ and $d=2$, a measure on the sphere of radius $t$ for $d=3$, and a distribution with compact support for $d\geq 4$ (see e.g. Theorem 5.28 of \cite{folland95}).
The definitions of $G^w$ and $G^h$ are not important for the present article. What is important is the form of their Fourier transforms: for any $t>0$ and $\xi \in \bR^d$,
\begin{equation}
\label{Fourier-G}
\cF G^w(t,\cdot)(\xi)=\frac{\sin(t|\xi|)}{|\xi|} \quad \mbox{and} \quad \cF G^h(t,\cdot)(\xi)=\exp\left(-\frac{t|\xi|^2}{2}\right).
\end{equation}

Due to the constant initial condition, the law of the solution $u(t,x)$ does not depend on $x$. In particular, the second moment $h(t)=E|u(t,x)|^2$ is independent of $x$, a fact which can also be observed from relation \eqref{second-moment} below. A natural question is to find out what is the Laplace transform of the function $h$. This calculation lies at the origin of our investigations, and was inspired by the beautiful proof of Theorem 1.8 of \cite{FK13}. A key observation made along this calculation is the fact that the Laplace transforms of $\cF G^w(t,\cdot)(\xi)$ and $\cF G^h(t,\cdot)(\xi)$ (as functions of $t$) have very similar forms. We exploited this fact in two different ways, which lead to the two main theorems below.


In case (i), we consider the functional
$$\cE(f) = \sup_{g\in \cF_d}\bigg\{\int_{\bR^d}
f(x)g^2(x)dx
-{1\over 2}\int_{\bR^d}\vert\nabla g(x)\vert^2dx\bigg\},$$
where
$${\cal F}_d=\left\{ g\in H^1(\bR^d); \int_{\bR^d} g^2(x)dx=1 \text{ and } \int_{\bR^d}|\nabla g(x)|^2 dx<\infty \right\}.$$
$H^1(\bR^d)$ is the Sobolev space of order $1$ on $\bR^d$, $\nabla g(x)=(\frac{\partial g}{\partial x_1}, \ldots,\frac{\partial g}{\partial x_d})$, and $\frac{\partial g}{\partial x_j}$ is the weak partial derivative of $g$ with respect to $x_j$ for $j=1,\ldots,d$. In Appendix \ref{appB}, we show that
\begin{equation}
\label{relation-2E}
\cE_2(f)=\cE(f(\sqrt{2}\,\cdot\,)),
\end{equation}
where $\cE_2(f)$ is the functional considered in \cite{huang-le-nualart15}:
$$\cE_2(f) = \sup_{g\in \cF_{2d}}\bigg\{\int_{\bR^d} \int_{\bR^d}
f(x_1-x_2)g^2(x_1,x_2)dx_1 dx_2
-{1\over 2}\int_{\bR^d} \int_{\bR^d}\vert\nabla g(x_1,x_2)\vert^2dx_1 dx_2\bigg\}.$$

In case (ii), we consider the functionals:
\begin{align*}
\cE(\Gamma)&=\sup_{h \in \cA_1}\bigg\{\int_{\bR}
 (h*h)(\xi)\mu(d\xi)-\frac{1}{2}\int_{\bR}|\xi|^2 |h(\xi)|^2 d\xi  \bigg\},\\
\cE_2(\Gamma)&=\sup_{h \in \cA_2}\bigg\{\int_{\bR}
 (h*h)(\xi,-\xi)\mu(d\xi)-\frac{1}{2}\int_{\bR}\int_{\bR}(|\xi_1|^2+|\xi_2|^2) |h(\xi_1,\xi_2)|^2 d\xi_1 d\xi_2  \bigg\}.
\end{align*}
where
$${\cal A}_d=\left\{ h:\bR^d \to \bC; \int_{\bR^d} |h(\xi)|^2d\xi=1, \int_{\bR^d}|\xi|^2|h(\xi)|^2 d\xi<\infty, \overline{h(\xi)}=h(-\xi) \forall \xi \in \bR^d\right\}.$$

Here is our first main result.

\begin{theorem}
\label{main-th1}
In case (i), suppose that there exists $\alpha>0$ such that
\begin{equation}
\label{scaling-f}
f(cx)=c^{-\alpha}f(x) \quad \mbox{for any}  \ c>0,x \in \bR^d,
\end{equation}
The upper second order Lyapunov exponent of the solution to equation \eqref{HAM} is:
$$\limsup_{t\to \infty}\frac{1}{t}\log E|u^w(t,x)|^2=
\left\{
\begin{array}{ll}
2^{\frac{1-\alpha}{3-\alpha}}
\Big(\cE_2(f)\Big)^{\frac{2-\alpha}{6-2\alpha}}=
2^{\frac{2-3\alpha}{6-2\alpha}}\Big(\cE(f)\Big)^{\frac{2-\alpha}{6-2\alpha}}
& \mbox{in case (i)} \\
2^{\frac{2H-1}{2H+1}}
\Big(\cE_2(\Gamma)\Big)^{\frac{H}{2H+1}}=2^{\frac{3H-2}{2H+1}}
\Big(\cE(\Gamma)\Big)^{\frac{H}{2H+1}} & \mbox{in case (ii).}
\end{array} \right.
$$
\end{theorem}

\begin{remark}
{\rm Using relation \eqref{Fourier-f-mu} between $f$ and $\mu$, we see that \eqref{scaling-f} is equivalent to the following scaling property of $\mu$: for any $c>0$ and for any Borel set $A\subset \bR^d$,
$$\mu(cA)=c^{\alpha}\mu(A).$$
We show in Appendix \ref{appA} that a tempered measure $\mu$ with this scaling property satisfies Dalang's condition \eqref{Dalang-cond} if and only if $\alpha<2$. Interestingly, the L\'evy measure $\nu$ of an $\alpha$-stable distribution in $\bR^d$ satisfies the scaling property $\nu(cA)=c^{-\alpha}\nu(A)$ for any $c>0$ and for any Borel set $A \subset \bR^d$ (see e.g. Theorem 14.3 of
of \cite{sato99}).}
\end{remark}

For our second result, we assume that $f$ is the Riesz kernel of order $\alpha$. We define
\begin{equation}
\label{def-a}
a=\left\{
\begin{array}{ll}
3-\alpha & \mbox{for equation \eqref{HAM}} \\
1-\alpha/2 & \mbox{for equation \eqref{PAM}}
\end{array} \right.
\end{equation}
\begin{equation}
\label{def-gamma}
\gamma=\left\{
\begin{array}{ll}
\log(2^{1-\alpha}\rho) & \mbox{for equation \eqref{HAM}} \\
\log \rho & \mbox{for equation \eqref{PAM}}
\end{array} \right.
\end{equation}
where $\rho$ is the constant introduced in \cite{bass-chen-rosen09}:
\begin{equation}
\label{def-rho}
\rho=\sup_{\|g\|_2=1}\left\{\int_{\bR^d}\int_{\bR^d}\frac{g(\xi)}{\sqrt{1+|\xi|^2}}
\frac{g(\eta)}{\sqrt{1+|\eta|^2}}\,C_{d,\alpha}|\xi-\eta|^{-(d-\alpha)}d\xi d\eta\right\},
\end{equation}
 {with $C_{d,\alpha}$ given by \eqref{def-C-d-alpha}.} Note that $\rho<\infty$ since $\alpha<2$ (see \cite{bass-chen-rosen09}).

 {Our second main result is the following.}

\begin{theorem}
\label{main-th2}
If $f(x)=|x|^{-\alpha}$ with $0<\alpha<d \wedge 2$, then
the second order Lyapunov exponent for the solution to equation \eqref{HAM} and for the solution to equation \eqref{PAM} is given by:
\begin{equation}
\label{lambda2}
\lambda_2=e^{\gamma/a}=\left\{
\begin{array}{ll} (2^{1-\alpha} \rho)^{1/(3-\alpha)} & \mbox{for equation \eqref{HAM}} \\
\rho^{2/(2-\alpha)} & \mbox{for equation \eqref{PAM}}
\end{array} \right.
\end{equation}
with constants $a$, $\gamma$ and $\rho$ defined
 by \eqref{def-a}, \eqref{def-gamma} and \eqref{def-rho}, respectively.
\end{theorem}

\begin{remark}(Comparison of the two results)
\label{comparison-remark}
 {
{\rm
When $f(x)=|x|^{-\alpha}=:R_{\alpha}(x)$, Theorems \ref{main-th1} and \ref{main-th2} give the same value for $\lambda_2$. To see this, on one hand, we express $\cE(R_{\alpha})$ in Theorem \ref{main-th1} using Lemma \ref{E-A-functional} (Appendix \ref{appB}) as follows:
\begin{equation}
\label{E-Riesz1}
\cE(R_{\alpha})=\cE(R_{\alpha},1/2)=2^{-\frac{\alpha}{\alpha-2}}\cE(R_{\alpha},1),
\end{equation}
where $\cE(f,A)$ is given by \eqref{def-E-A}.
On the other hand, by Theorem 1.5 of \cite{bass-chen-rosen09}, we know that:
\begin{equation}
\label{E-Riesz2}
\rho=\Lambda_{\alpha}^{1-\alpha/2}=\cE(R_{\alpha},1)^{1-\alpha/2},
\end{equation}
(Note that there is a small error in Theorem 1.5 of \cite{bass-chen-rosen09} which states that $\rho=(2\pi)^{-d}\Lambda_{\alpha}^{1-\alpha/2}$. This error is due to the fact that in the first line of equation (7.22) of \cite{bass-chen-rosen09} one should have $(2\pi)^{-dp}$ instead of $(2\pi)^{-d(p+1)}$, since
$\cF (f^p)=(2\pi)^{-d(p-1)}(\cF f)^{*p}$.) A direct calculation based on \eqref{E-Riesz1} and \eqref{E-Riesz2} shows that:
$$\left(2^{1-\alpha}\rho\right)^{1/(3-\alpha)}=
2^{\frac{2-3\alpha}{6-2\alpha}}\Big(\cE(R_{\alpha})\Big)^{\frac{2-\alpha}{6-2\alpha}}.$$
}
}
\end{remark}

\begin{remark}(Equations with fractional power of Laplacian)
{\rm Theorems \ref{main-th1} and \ref{main-th2} can be extended (with the same proof) to equations \eqref{HAM} and \eqref{PAM} in which the Laplacian $\Delta$ is replaced by $-(-\Delta)^{\beta/2}$ for $\beta \in (0,2]$, provided the measure $\mu$ satisfies
$$\int_{\bR^d}\frac{1}{1+|\xi|^{\beta}}\mu(d\xi)<\infty,$$
which is equivalent to $\alpha<\beta$. In this case, formula \eqref{Fourier-G} remains valid with $|\xi|$ replaced by $|\xi|^{\beta/2}$ (see e.g. Section 3.2 of \cite{FKN11} for the wave equation), and existence of the mild solution can be proved similarly to \cite{balan-song16}. In case (i), the statement of Theorem \ref{main-th1} becomes:
$$\limsup_{t\to \infty}\frac{1}{t}\log E|u^w(t,x)|^2=
2^{\frac{1}{3\beta-2\alpha}[\beta(1-\alpha)+\frac{\alpha(\beta-\alpha)}{\alpha-2}]}
\Big(\cE(f)\Big)^{\frac{\beta-\alpha}{3\beta-2\alpha}}.
$$
A similar relation holds for case (ii) with $\alpha$ replaced by $2-2H$ and $\cE(f)$ replaced by $\cE(\Gamma)$. Relation \eqref{lambda2} holds with $\gamma$ given by \eqref{def-gamma}, $a$ replaced by:
\begin{equation}
a=\left\{
\begin{array}{ll}
3-2\alpha/\beta & \mbox{for equation \eqref{HAM}} \\
1-\alpha/\beta & \mbox{for equation \eqref{PAM}}
\end{array} \right.
\end{equation}
and $\rho$ replaced by
\begin{equation}
\label{def-rho'}
\rho=\sup_{\|g\|_2=1}\left\{\int_{\bR^d}\int_{\bR^d}\frac{g(\xi)}{\sqrt{1+|\xi|^\beta}}
\frac{g(\eta)}{\sqrt{1+|\eta|^\beta}}\,C_{d,\alpha}|\xi-\eta|^{-(d-\alpha)}d\xi d\eta\right\}.
\end{equation}
We refer the reader to \cite{Song17}  {for a study of a general} parabolic Anderson models with space-time colored noise, and to \cite{CHSS}  {for the precise calculation of the Lyapunov exponents of order $p\geq 2$} of the solution to a fractional parabolic Anderson model.
}
\end{remark}

\begin{remark}(Equations with space-time white noise)
{\rm With minor modifications, the proof of Theorem \ref{main-th2} can also be used for the case of equations \eqref{HAM} and \eqref{PAM} with space-time white noise in spatial dimension $d=1$. In this case,  $\cH=L^2(\bR_{+} \times \bR)$,
$$\langle \varphi,\psi \rangle_{\cH}=\int_{0}^{\infty}\int_{\bR}\varphi(t,x)\psi(t,x)dxdt,$$
and  relation \eqref{inner-product-H} holds with
$\mu(d\xi)=(2\pi)^{-1} d\xi$, by Plancherel theorem. Intuitively, this corresponds to the case of Riesz kernel with critical exponent $\alpha=1$. The results in Section \ref{section-proof2} below hold (with the same arguments) with constants:
$$a=\left\{
\begin{array}{ll}
2 & \mbox{for equation \eqref{HAM}} \\
1/2 & \mbox{for equation \eqref{PAM}}
\end{array} \right.
\quad \mbox{and} \quad
\gamma=\rho=\log(1/2).
$$
The proof of Lemma \ref{limit-Jn-tau} below shows that in this case,
$E[J_n^w(\tau)]=E[J_n^h(\tau)]=T_n$, where
$$T_n=\frac{1}{(2\pi)^{n}}\int_{\bR^n}\frac{1}{1+|\xi_1|^2} \frac{1}{1+|\xi_1+\xi_2|^2}\ldots \frac{1}{1+|\xi_1+\ldots+\xi_n|^2} d\xi_1 \ldots d\xi_n=\frac{1}{2^n}.$$
We obtain the well-known results:
$$\lambda_2=e^{\gamma/a}=(1/2)^{1/a}=\left\{
\begin{array}{ll}
1/\sqrt{2} & \mbox{for equation \eqref{HAM}} \\
1/4 & \mbox{for equation \eqref{PAM}}
\end{array} \right. .
$$
}
\end{remark}

 {
This article is organized as follows. In Section \ref{section-th1}, we compute the Laplace transforms of the second moment of the solutions to equations \eqref{HAM} and \eqref{PAM}. In Section \ref{section-connection}, we derive the connection between the Laplace transforms in the hyperbolic and parabolic case, leading to the proof of Theorem \ref{main-th1}. In Section \ref{section-proof2} we give the proof of Theorem \ref{main-th2}. For this, we use the special form (when $f$ is the Riesz kernel) of the Laplace transform of the $n$-th term $J_n(t)$ appearing in the series representation of the second moment of the solution, and a key result borrowed from \cite{bass-chen-rosen09}. Finally, we present some auxiliary results in Appendix \ref{appA}, while in Appendix \ref{appB}, we derive some scaling properties of the functionals $\cE$ and $\cE_2$ which are needed in the sequel.}

\section{Laplace transforms}
\label{section-th1}

In this section, we gather some useful facts about the Laplace transforms of various deterministic functions which are used in this article. We also recall some basic facts about the existence of the solution and the calculation of its second moment. The scaling property \eqref{scaling-f} is not needed for the results presented in this section.

Using the same method as in \cite{balan-song16} (for case (i)) and in \cite{BJQ} (for case (ii)), it can be proved that the solution $u$ to either one of equations \eqref{HAM} or \eqref{PAM} has the Wiener chaos expansion:
\begin{equation}
\label{chaos-exp}
u(t,x)=1+\sum_{n\geq 1}I_n(f_n(\cdot,t,x)),
\end{equation}
where $I_n$ is the multiple integral of order $n$ with respect to $W$, and
$$f_n(t_1,x_1,\ldots,t_n,x_n,t,x)=G(t-t_n,x-x_n) \ldots G(t_2-t_1,x_2-x_1)1_{\{0<t_1<\ldots<t_n<t\}}$$
with $G=G^w$ for equation \eqref{HAM} with $d \leq 2$, respectively $G=G^h$ for equation \eqref{PAM}. Relation \eqref{chaos-exp} holds also for the solution to equation \eqref{HAM} with $d \geq 3$, but in this case $f_n(t_1,\cdot,\ldots,t_n,\cdot,t,x)$ is a distribution in $\cS'(\bR^{nd})$.
In all cases, the Fourier transform of $f_n(t_1,\cdot,\ldots,t_n,\cdot,t,x)$
is
\begin{eqnarray}
\label{Fourier-fn}
\lefteqn{
\cF f_n(t_1,\cdot, \ldots, t_n,\cdot,t,x)(\xi_1, \ldots,\xi_n)=e^{-i (\xi_1 +\ldots+\xi_n)\cdot x}\cF G(t_2-t_1,\cdot)(\xi_1) } \\
\nonumber
& & \cF G(t_3-t_2,\cdot)(\xi_1+\xi_2)\ldots \cF G(t-t_n,\cdot)(\xi_1+\ldots+\xi_n)1_{\{0<t_1<\ldots<t_n<t\}}.
\end{eqnarray}

By the orthogonality of the Wiener chaos spaces,
\begin{equation}
\label{second-moment}
E|u(t,x)|^2=1+\sum_{n \geq 1} E|I_n(f_n(\cdot,t,x))|^2=1+\sum_{n \geq 1} J_n(t),
\end{equation}
where $J_n(t)=n!\,\|\widetilde{f}_n(\cdot,t,x)\|_{\cH^{\otimes n}}^{2}$ and $\widetilde{f}_n(\cdot,t,x)$ is the symmetrization of $f_n(\cdot,t,x)$. We let $J_n(0)=1$.

\begin{lemma}
\label{norm-fn}
For any $t>0$, $x \in \bR^d$ and $n \geq 1$,
\begin{eqnarray}
\label{def-Jn}
J_n(t)&=&
\int_{0<t_1<\ldots<t_n<t} \int_{\bR^{nd}}|\cF G(t_2-t_1,\cdot)(\xi_1)|^2 |\cF G(t_3-t_2,\cdot)(\xi_1+\xi_2)|^2 \ldots  \\
\nonumber
& & \quad \quad \quad |\cF G(t-t_n,\cdot)(\xi_1+\ldots+\xi_n)|^2 \mu(d\xi_1) \ldots \mu(d\xi_n)dt_1 \ldots dt_n,
\end{eqnarray}
and $G=G^{w}$ for equation \eqref{HAM}, respectively $G=G^h$ for equation \eqref{PAM}.
We denote by $J_n^w(t)$ and $J_n^h(t)$ the integral given by \eqref{def-Jn} with $G=G^{w}$, respectively $G=G^h$.
\end{lemma}

\noindent {\bf Proof:} The norm in the space $\cH^{\otimes n}$ is computed using Fourier transforms, as in \eqref{inner-product-H}. Note that $\cF \widetilde{f}_n(t_1,\cdot, \ldots, t_n,\cdot,t,x)(\xi_1, \ldots,\xi_n)$ is the sum over all permutations $\rho$ on $\{1,\ldots,n\}$ of $n!$ terms which have a similar expression to \eqref{Fourier-fn} in which $(t_1,\ldots,t_n)$ and $(\xi_1,\ldots,\xi_n)$ are replaced by $(t_{\rho(1)},\ldots,t_{\rho(n)})$, respectively $(\xi_{\rho(1)},\ldots,\xi_{\rho(n)})$.
Due to the presence of the indicator of $\{0<t_{\rho(1)}<\ldots<t_{\rho(n)}<t\}$ in all such terms, in the computation of the squared modulus of $\cF \widetilde{f}_n(t_1,\cdot, \ldots, t_n,\cdot,t,x)(\xi_1, \ldots,\xi_n)$, the mixed terms corresponding to different permutations $\rho$ and $\sigma$ vanish, and we are left only with the $n!$ terms corresponding to permutations $\rho=\sigma$. The conclusion follows recalling that the definition of the symmetrization $\widetilde{f}_n(\cdot,t,x)$ contains the factor $1/n!$ $\Box$

\vspace{3mm}


For any $\xi \in \bR^d$ and $\beta>0$, we consider the Laplace transforms:
$$\cI_{\beta}^w(\xi)=\int_{0}^{\infty}e^{-\beta t}\, |\cF G^w(t,\cdot)(\xi)|^2 \,dt \quad \mbox{and} \quad
\cI_{\beta}^h(\xi)=\int_{0}^{\infty}e^{-\beta t}\, |\cF G^h(t,\cdot)(\xi)|^2 \,dt.$$

\begin{lemma}
\label{lemma-I-beta}
For any $\beta>0$ and for any $\xi \in \bR^d$, we have:
$$\cI_{\beta}^w(\xi)=\frac{1}{2\beta}\cdot\frac{1}{\frac{\beta^2}{4}+|\xi|^2} \quad \mbox{and} \quad
\cI_{\beta}^h(\xi)=\frac{1}{\beta+|\xi|^2}.$$
\end{lemma}

\noindent {\bf Proof:} The result for $\cI^h(\beta)$ is clear.
For $\cI^{w}(\beta)$, we use the identity:
$$\int_0^{\infty}e^{-ux}\sin^2 xdx=\frac{2}{u(u^2+4)}, \quad u>0,$$
which can be deduced from $\int_0^{\infty}e^{-ux} \cos x dx=\frac{u}{u^2+1}$, since  $2\sin^2 x=1-\cos (2x)$.
Let $\xi \in \bR$ be arbitrary. Using the change of variable $x=t|\xi|$, we have:
\begin{eqnarray*}
\cI_{\beta}^w(\xi)&=&\frac{1}{|\xi|^2}\int_0^{\infty}e^{-\beta t} \sin^2(t|\xi|)dt=\frac{1}{|\xi|^3}\int_0^{\infty}e^{-\beta x/|\xi|}\sin^2 x dx\\
&=& \frac{1}{|\xi|^3} \cdot
\frac{2}{\frac{\beta}{|\xi|}(\frac{\beta^2}{|\xi|^2}+4)}=\frac{2}{\beta(\beta^2+4|\xi|^2)}=
\frac{1}{2\beta}\cdot\frac{1}{\frac{\beta^2}{4}+|\xi|^2}
\end{eqnarray*}
$\Box$

We now compute the Laplace transform of the second moment of the solution:
$$L(\beta)=\int_{0}^{\infty}e^{-\beta t}E|u(t,x)|^2 dt, \quad \beta>0.$$

\begin{lemma}
\label{Laplace-second-mom}
For any $\beta>0$ and $n\geq 1$, we have:
\begin{equation}
\label{integral-Jn}
\int_{0}^{\infty}e^{-\beta t}J_n(t)dt= \frac{1}{\beta}\int_{\bR^{nd}}\cI_{\beta}(\xi_1)\cI_{\beta}(\xi_1+\xi_2)
\ldots \cI_{\beta}(\xi_1+\ldots+\xi_n)\mu(d\xi_1)\ldots \mu(d\xi_n).
\end{equation}
Therefore, for any $\beta>0$,
$$L(\beta)= \frac{1}{\beta}\Big(1+\sum_{n\geq 1}\int_{\bR^{nd}}\cI_{\beta}(\xi_1)\cI_{\beta}(\xi_1+\xi_2)
\ldots \cI_{\beta}(\xi_1+\ldots+\xi_n)\mu(d\xi_1)\ldots \mu(d\xi_n)\Big).$$
\end{lemma}

\noindent {\bf Proof:} For the first equality, we use \eqref{def-Jn} and then write $e^{-\beta
t}=e^{-\beta t_1} e^{-\beta (t_2-t_1)} \ldots e^{-\beta(t-t_n)}.$
Using the change of variable $u=t_1, u_1=t_2-t_1, \ldots,
u_n=t-t_n$, and Fubini's theorem, we see that
\begin{eqnarray*}
\int_{0}^{\infty}e^{-\beta t}J_n(t)dt&=&
\int_{\bR^{nd}}\left(\int_{0}^{\infty}e^{-\beta u}du\right) \left(\int_0^{\infty}e^{-\beta u_1}|\cF G(u_1,\cdot)(\xi_1)|^2 du_1 \right)\ldots \\
& & \left(\int_0^{\infty}e^{-\beta u_n}|\cF
G(u_n,\cdot)(\xi_1+\ldots+\xi_n)|^2 du_n \right)\mu(d\xi_1)\ldots
\mu(d\xi_n).
\end{eqnarray*}
The second equality follows by \eqref{second-moment} and Fubini's theorem:
$$\int_{0}^{\infty}e^{-\beta t}E|u(t,x)|^2 dt=\int_{0}^{\infty}e^{-\beta t}dt+\sum_{n \geq 1}\int_{0}^{\infty}e^{-\beta t}J_n(t)dt.$$
$\Box$

\section{Proof of Theorem \ref{main-th1}}
\label{section-connection}

In this section, we prove Theorem \ref{main-th1}. More precisely, we will derive the second order upper Lyapunov exponent of $u^w(t,x)$, using the {\em global} asymptotic behaviour of the second moment of solutions $\{u^{h,\beta}\}_{\beta>0}$ to a class of parabolic Anderson models perturbed by a family $\{W^{\beta}\}_{\beta>0}$ of Gaussian noise processes,
whose precise definition is given below.

The key observation is the fact that the Laplace transforms $E|u^w(t,x)|^2$ and $E|u^h(t,x)|^2$ have {\em exactly} the same form, relying on the respective quantities $\cI_{\beta}^{w}(\xi)$ and $\cI_{\beta}^h(\xi)$ (see Lemma \ref{Laplace-second-mom}). We denote these Laplace transforms by $L^w$, respectively $L^{h}$.

Note that by Lemma \ref{lemma-I-beta}, for any $\beta>0$ and for any $\xi \in \bR^d$, we have
\begin{equation}
\label{I-wh}
\cI_\beta^w(\xi)=\frac{1}{2\beta} \cI_{\gamma(\beta)}^h(\xi),
\end{equation}
where $\gamma(\beta):=\beta^2/4$. We denote $\mu_{\beta}(d\xi)=(2\beta)^{-1}\mu(d\xi)$.

The Fourier transform in $\cS'(\bR^d)$ of the measure $\mu_{\beta}$ is the function $f_{\beta}=(2\beta)^{-1}f$ in case (i), respectively the distribution $\Gamma_{\beta}=(2\beta)^{-1}\Gamma$ in case (ii).

Let $W^{\beta}=\{W^{\beta}(\varphi); \varphi \in \cH^{\beta}\}$ be an isonormal Gaussian process with covariance $$E[W^{\beta}(\varphi)W^{\beta}(\psi)]=\langle \varphi,\psi \rangle_{\cH_{\beta}},$$ where $\cH_{\beta}$ is the completion of $C_0^{\infty}(\bR_{+} \times \bR^d)$ with respect to the inner product $\langle \cdot,\cdot \rangle_{\cH_{\beta}}$ defined by \eqref{inner-product-H} with $\mu$ replaced by $\mu_{\beta}$.

Let $u^{h,\beta}$ be the solution to equation \eqref{PAM} with $W$ replaced by $W^{\beta}$. We denote by $L^{h,\beta}$ the Laplace transform of the second moment of the solution $u^{h,\beta}$:
$$L^{h,\beta}(\lambda)=\int_{0}^{\infty}e^{-\lambda t}E|u^{h,\beta}(t,x)|^2 dt, \quad \lambda>0.$$

\begin{lemma}
For any $\beta>0$ and for any $x \in \bR^d$, we have:
\begin{equation}
\label{connection-xy}
\int_0^{\infty} e^{-\beta t} E|u^{w}(t,x)|^2 dt=\frac{\beta}{4} \int_0^{\infty} e^{- t\beta^2/4} E|u^{h,\beta}(t,x)|^2 dt.
\end{equation}
\end{lemma}

\noindent {\bf Proof:} By Lemma \ref{Laplace-second-mom} and relation \eqref{I-wh}, the
Laplace tranform $L^w(\beta)$ is equal to:
\begin{align*}
& \frac{\gamma(\beta)}{\beta} \frac{1}{\gamma(\beta)}\Big(1+\sum_{n\ge1}\int_{\bR^{nd}}
\cI^h_{\gamma(\beta)}(\xi_1)\cI^h_{\gamma(\beta)}(\xi_1+\xi_2) \ldots \cI^h_{\gamma(\beta)}(\xi_1+\ldots+\xi_n) \mu_{\beta}(d\xi_1)\ldots \mu_{\beta}(d\xi_n)\Big) \notag\\
&=\frac{\beta}{4} L^{h,\beta}(\gamma(\beta)).
\end{align*}
Note that the last equality is due also to Lemma \ref{Laplace-second-mom} applied to the solution $u^{h,\beta}$. $\Box$

\vspace{3mm}

The following result shows how to derive the asymptotic behaviour of $E|u^w(t,x)|^2$, assuming that we can control the behaviour of $E|u^{h,\beta}(t,x)|^2$ for all $\beta>0$.

\begin{theorem}
\label{connection}
Assume that for any $\beta>0$,
\begin{equation}
\label{def-A-beta}
\lambda(\beta):=\limsup_{t\to \infty}\frac{1}{t} \log E|u^{h,\beta}(t,x)|^2<\infty.
\end{equation}
Suppose that $\lambda:(0,\infty) \to [0,\infty)$ is continuous and strictly decreasing and satisfies $\lim_{\beta \to 0+}\lambda(\beta) \in (0,\infty]$ and $\lim_{\beta \to \infty}\lambda(\beta)=c \in [0,\infty)$. Then
\begin{equation}
\label{swelim}
\limsup_{t\to \infty}\frac1t \log E|u^{w}(t,x)|^2=\beta_0,
\end{equation}
where $\beta_0$ is the unique solution of the equation $4\lambda(\beta)=\beta^2$ in $(0,\infty)$.
\end{theorem}

\noindent{\bf Proof:} We consider the non-decreasing function $h_{\beta}(t)=E|u^{h,\beta}(t,x)|^2$. By Lemma \ref{lemma-loglim} (Appendix \ref{appA}),
$$\int_0^{\infty} e^{- t\beta^2/4} E|u^{h,\beta}(t,x)|^2 dt<\infty \quad \mbox{if} \quad \frac{\beta^2}{4}>\lambda(\beta)$$
and
$$\int_0^{\infty} e^{- t\beta^2/4} E|u^{h,\beta}(t,x)|^2 dt=\infty \quad \mbox{if} \quad \frac{\beta^2}{4}<\lambda(\beta).$$

\noindent Define $g(\beta)=4\lambda(\beta)-\beta^2$ for any $\beta>0$. Then $g$ is a continuous strictly decreasing function on $(0, \infty)$, which satisfies $\lim_{\beta \to 0+}g(\beta)=4\lim_{\beta \to 0+}\lambda(\beta)\in (0,\infty]$ and $\lim_{\beta\to \infty} g(\beta)=4c-\infty=-\infty$. Hence, the equation $g(\beta)=0$ has a unique solution $\beta_0$ in $(0,\infty)$. Moreover, $g(\beta)>0$ for all $\beta<\beta_0$ and $g(\beta)<0$ for all $\beta>\beta_0$. This means that $\frac{\beta^2}{4}<\lambda(\beta)$ for all $\beta<\beta_0$ and $\frac{\beta^2}{4}>\lambda(\beta)$ for all $\beta>\beta_0$.
We obtain that:
$$\int_0^{\infty} e^{- \beta t} E|u^{w}(t,x)|^2 dt<\infty \quad \mbox{if} \quad \beta>\beta_0$$
and
$$\int_0^{\infty} e^{- \beta t} E|u^{w}(t,x)|^2 dt=\infty \quad \mbox{if} \quad \beta<\beta_0.$$
We now apply again Lemma \ref{lemma-loglim} (Appendix \ref{appA}) to the non-decreasing function $h(t)=E|u^w(t,x)|^2$. Relation \eqref{swelim} follows. $\Box$

\vspace{3mm}

We are now ready to give the proof of Theorem \ref{main-th1}.

\vspace{3mm}

\noindent {\bf Proof of Theorem \ref{main-th1}:} We first treat case (i). By Theorem 1.2 of \cite{huang-le-nualart15}, we know that
$$\lambda(\beta)=\lim_{t \to \infty}\frac{1}{t}\log E|u^{h,\beta}(t,x)|^2=\cE_2(f_{\beta}).$$
Using the definition of $f_{\beta}$ and Lemma \ref{scaling-E} (Appendix \ref{appB}), we see that
$$\lambda(\beta)=(2\beta)^{-\frac{2}{2-\alpha}}\cE_2(f)
=2^{-\frac{2+\alpha}{2-\alpha}}\beta^{-\frac{2}{2-\alpha}}\cE(f).$$

The function $\lambda$ is continuous and strictly decreasing and satisfies $\lim_{\beta \to 0+}\lambda(\beta)=\infty$ and $\lim_{\beta \to \infty}\lambda(\beta)=0$. The unique solution of the equation $4 \lambda(\beta)=\beta^2$ in $(0,\infty)$ is
$$\beta_0=2^{\frac{1-\alpha}{3-\alpha}}\Big(\cE_2(f) \Big)^{\frac{2-\alpha}{6-2\alpha}}=2^{\frac{2-3\alpha}{6-3\alpha}}\Big(\cE(f)\Big)^{\frac{2-\alpha}{6-2\alpha}}.$$
The conclusion follows by Theorem \ref{connection}.

Next we consider case (ii). By Theorem 1.2 of \cite{huang-le-nualart15} and Lemma \ref{scaling-E-Gamma} (Appendix \ref{appB}),
$$\lambda(\beta)=\lim_{t \to \infty}\frac{1}{t}\log E|u^{h,\beta}(t,x)|^2=\cE_2(\Gamma_{\beta})=(2\beta)^{-1/H}\cE_2(\Gamma)=
2^{-\frac{2-H}{H}}\beta^{-1/H}\cE(\Gamma).$$
The conclusion follows as in case (i).
$\Box$

\section{Proof of Theorem \ref{main-th2}}
\label{section-proof2}

In this section, we give the proof of Theorem \ref{main-th2}.
The key observation is that,  {when $f$ is the Riesz kernel of order $\alpha$}, the integral appearing on the right-hand side of \eqref{integral-Jn} is related to the Riesz potential $\zeta(t)=\int_{0}^t |B_s|^{-\alpha}ds$ of a $d$-dimensional Brownian motion $B=(B_t)_{t \geq 0}$.
More precisely, if $\tau_1$ is an exponential random variable with mean $1$, independent of $B$, by Lemma 2.2 of \cite{bass-chen-rosen09},
$E[\zeta(\tau_1)^n]=n!\,T_n$,
where
$$T_n= C_{d,\alpha}^n \int_{\bR^{nd}} \frac{1}{1+|\xi_1|^2} \ldots \frac{1}{1+|\xi_1+\xi_2|^2}\ldots \frac{1}{1+|\xi_1+\ldots+\xi_n|^2} \prod_{i=1}^{n}|\xi_i|^{\alpha-d}d\xi_1 \ldots d \xi_n,$$
 {and the constant $C_{d,\alpha}$ is given by \eqref{def-C-d-alpha}}.
The exact asymptotic behaviour of $E[\zeta(\tau_1)^n]/n!=T_n$ is given by Theorem 2.2 of \cite{bass-chen-rosen09}:
\begin{equation}
\label{th2-2-bass}
\lim_{n \to \infty}\frac{1}{n} \log T_n=\log \rho,
\end{equation}
where $\rho$ is the constant given by \eqref{def-rho}. This leads to the following result.

\begin{lemma}
\label{limit-Jn-tau}
Let $\tau$ be an exponential random variable with mean $1$. Then
$$\lim_{n \to \infty}\frac{1}{n}\log E[J_n(\tau)]=\gamma,$$
where $\gamma$ is the constant given by \eqref{def-gamma}.
\end{lemma}

\noindent {\bf Proof:} Using \eqref{integral-Jn} and Lemma \ref{lemma-I-beta} with $\beta=1$, we obtain:
\begin{eqnarray*}
\lefteqn{E[J_n^w(\tau)]=\int_{0}^{\infty}e^{-t}J_n^w(t)dt} \\
& & =
C_{d,\alpha}^n \int_{\bR^{nd}}\cI_{1}^w(\xi_1)\cI_{1}^w(\xi_1+\xi_2)\ldots \cI_{1}^w(\xi_1+\ldots+\xi_n)\prod_{i=1}^{n}|\xi_{i}|^{\alpha-d}d\xi_1 \ldots d\xi_n\\
& & = C_{d,\alpha}^n \int_{\bR^{nd}} \frac{2}{1+4|\xi_1|^2} \ldots \frac{2}{1+4|\xi_1+\xi_2|^2}\ldots \frac{2}{1+4|\xi_1+\ldots+\xi_n|^2} \prod_{i=1}^{n}|\xi_i|^{\alpha-d}d\xi_1 \ldots d \xi_n\\
& & = 2^{n(1-\alpha)}T_n.
\end{eqnarray*}
A similar calculation shows that $E[J_n^h(\tau)]=T_n$. The conclusion follows by \eqref{th2-2-bass}. $\Box$

\vspace{3mm}

The next result shows that the terms $J_n(t)$ have a scaling property in $t$.

\begin{lemma}
\label{scaling}
Let $a$ be the constant given by \eqref{def-a}. For any $t>0$ and $n \geq 0$, we have:
\begin{equation}
\label{scaling-J}
J_n(t)=t^{an}J_n(1).
\end{equation}
\end{lemma}

\noindent {\bf Proof:} This follows by a change of variables. For the time variables,
we let $s_i=t_i/t$ for $i=1,\ldots,n$ for both equations. For the space variables, we let
$\eta_i=t\xi_i$ for $i=1,\ldots,n$ for equation \eqref{HAM}, and $\eta_i=\sqrt{t}\xi_i$ for $i=1,\ldots,n$ for equation \eqref{PAM}. $\Box$

\vspace{3mm}

From Lemma \ref{limit-Jn-tau} and the scaling property of $J_n(t)$, we deduce the asymptotic behaviour of $J_n(1)$.

\begin{lemma}
Let $a$ and $\gamma$ be the constants given by \eqref{def-a}, respectively \eqref{def-gamma}. Then
\begin{equation}
\label{limit-J}
\lim_{n \to \infty}\frac{1}{n}\log \Big(\Gamma(an+1)J_n(1) \Big)=\gamma.
\end{equation}
\end{lemma}

\noindent {\bf Proof:} Let $\tau$ be an exponential random variable with mean $1$. By Lemma \ref{scaling}, $J_n(\tau)=\tau^{an}J_n(1)$. Then $$E[J_n(\tau)]=E[\tau^{an}]J_n(1)=\Gamma(an+1)J_n(1).$$
The conclusion follows by Lemma \ref{limit-Jn-tau}.
$\Box$

\vspace{3mm}

\noindent {\bf Proof of Theorem \ref{main-th2}:} Let $\varepsilon>0$ be arbitrary. By \eqref{scaling-J} and \eqref{limit-J}, there exists an integer $N=N_{\varepsilon} \geq 1$ such that for any $n \geq N$ and $t >0$,
\begin{equation}
\label{Jn-bound}
\frac{e^{n(\gamma-\varepsilon)}t^{an}}{\Gamma(an+1)} \leq J_n(t) \leq \frac{e^{n(\gamma+\varepsilon)}t^{an}}{\Gamma(an+1)}.
\end{equation}

{\em Step 1.} (Upper bound) Let $c_1=e^{(\gamma+\varepsilon)/a}$. By the upper bound in \eqref{Jn-bound}, for any $t>0$,
\begin{eqnarray*}
\sum_{n \geq N}J_n(t) & \leq & \sum_{n \geq N}\frac{(c_1 t)^{an}}{\Gamma(an+1)} \leq \sum_{n \geq 0}\frac{(c_1 t)^{an}}{\Gamma(an+1)}=: A_t^{(1)}.
\end{eqnarray*}
By Lemma \ref{lemmaA} (Appendix \ref{appA}), for any $\delta>0$, there exists some $t_{\delta,1}>1$ such that
$$A_t^{(1)} \leq e^{t(c_1+\delta)} \quad \mbox{for all} \quad t\geq t_{\delta,1}.$$
We denote $C_1=\sum_{n<N}J_n(1)$. We can find a value $t_{\delta,1}'>t_{\delta,1}$ large enough, such that
$$\sum_{n \geq 0}J_n(t)=\sum_{n<N}t^{an}J_n(1)+\sum_{n \geq N}J_n(t) \leq C_1 t^{aN}+e^{t(c_1+\delta)} \leq 2 e^{t(c_1+\delta)}$$
for all $t \geq t_{\delta,1}'$. Hence,
$$\frac{1}{t} \log \sum_{n \geq 0}J_n(t) \leq \frac{1}{t}\log 2+ c_1+\delta\leq c_1+2\delta$$
for all $t \geq t_{\delta,1}''$, for a value $t_{\delta,1}''>t_{\delta,1}'$ large enough. This implies that
\begin{equation}
\label{UB}
\limsup_{t \to \infty}\frac{1}{t}\log \sum_{n \geq 0}J_n(t) \leq c_1.
\end{equation}

{\em Step 2.} (Lower bound) Let $c_2=e^{(\gamma-\varepsilon)/a}$. By the lower bound in \eqref{Jn-bound}, for any $t>1$,
$$\sum_{n \geq N}J_n(t) \geq \sum_{n \geq N}\frac{(c_2 t)^{an}}{\Gamma(an+1)}=A_t^{(2)}-\sum_{n<N}\frac{(c_2 t)^{an}}{\Gamma(an+1)} \geq A_t^{(2)}-C_2 t^{aN},$$
where $$A_t^{(2)}=\sum_{n \geq 0} \frac{(c_2 t)^{an}}{\Gamma(an+1)} \quad \mbox{and} \quad C_2=\sum_{n<N}\frac{c_2^{an}}{\Gamma(an+1)}.$$
By Lemma \ref{lemmaA} (Appendix \ref{appA}), for any $\delta>0$ there exists some $t_{\delta,2}>1$ such that
$$A_{t}^{(2)} \geq e^{t(c_2-\delta)} \quad \mbox{for all} \quad t \geq t_{\delta,2}.$$
Therefore, we can find a value $t_{\delta,2}'>t_{\delta,2}$ such that for any $t \geq t_{\delta,2}'$,
$$\sum_{n \geq 0}J_n(t) \geq \sum_{n \geq N}J_n(t) \geq e^{t(c_2-\delta)}-C_2 t^{aN} \geq \frac{1}{2}e^{t(c_2-\delta)}.$$
It follows that
$$\frac{1}{t}\log \sum_{n \geq 0}J_n(t) \geq -\frac{\log 2}{t}+c_2-\delta \geq c_2- 2 \delta$$
for any $t \geq t_{\delta,2}''$, for a value $t_{\delta,2}''>t_{\delta,2}'$ large enough. This implies that
\begin{equation}
\label{LB}
\liminf_{t \to \infty}\frac{1}{t} \log \sum_{n \geq 0}J_n(t) \geq c_2.
\end{equation}

{\em Step 3.} (Conclusion) From \eqref{UB} and \eqref{LB}, we obtain that for any $\varepsilon>0$,
$$e^{(\gamma-\varepsilon)/a} \leq \liminf_{t \to \infty}\frac{1}{t} \log \sum_{n \geq 0}J_n(t) \leq \limsup_{t \to \infty}\frac{1}{t} \log \sum_{n \geq 0}J_n(t) \leq e^{(\gamma+\varepsilon)/a}.$$
The conclusion follows taking $\varepsilon \downarrow 0$. $\Box$

\appendix
\section{Some  {auxiliary} results}
\label{appA}

\begin{lemma}
\label{scaling-mu}
Let $\mu$ be a tempered measure on $\bR^d$ such that $\mu(cA)=c^{\alpha}\mu(A)$ for all $c>0$ and for any Borel set $A \subset \bR^d$. Then $\mu$ satisfies \eqref{Dalang-cond} if and only if $\alpha<2$.
\end{lemma}

\noindent {\bf Proof:} Since $\mu$ is locally integrable, \eqref{Dalang-cond} is equivalent to $I:=\int_{|\xi|>1}|\xi|^{-2}\mu(d\xi)<\infty$. Denote
$B_k:=\{\xi \in \bR^d: |\xi|\leq k\}$ and $A_k=B_{k+1}\backslash B_{k}$ for $k=1, 2,\dots$. Then
$$\sum_{k \geq 1}\frac{1}{(k+1)^2}\mu(A_k) \leq I \leq \sum_{k \geq 1}\frac{1}{k^2}\mu(A_k).$$
By the scaling property of $\mu$, $\mu(B_k)=k^{\alpha}\mu(B_1)$ and  $\mu(A_k)=[(k+1)^{\alpha}-k^{\alpha}]\mu(B_1)$.
Then
$$\mu(B_1)\sum_{k \geq 1} \frac{(k+1)^{\alpha}-k^{\alpha}}{(k+1)^2} \leq I \leq \mu(B_1)\sum_{k \geq 1} \frac{(k+1)^{\alpha}-k^{\alpha}}{k^2}.$$
Note that $\alpha x^{\alpha-1}\leq (x+1)^{\alpha}-x^{\alpha}\leq C(x^{\alpha-1}+1)$ for all $x>0$, where $C>0$ is a constant depending on $\alpha$. Hence
$I<\infty$ if and only if $\sum_{k \geq 1}k^{\alpha-3}<\infty$, i.e. $\alpha<2$. $\Box$

\begin{lemma}
\label{lemma-loglim}
Let $h:[0,\infty) \to [0,\infty)$ be a function such that  $\int_{0}^{\infty}e^{-\eta t}h(t)dt<\infty$ for some $\eta>0$.

a) Then
\begin{equation}
\label{def-A}
A=\inf\{\eta>0; \int_{0}^{\infty}e^{-\eta t}h(t)dt<\infty\}
\end{equation} if and only if
\begin{equation}
\label{finite-int1}
\int_0^\infty e^{-\eta t} h(t)dt<\infty \ for \ all \ \eta>A
\end{equation}
and
\begin{equation}
\label{finite-int2}
\int_0^\infty e^{-\eta t} h(t)dt=\infty \ for \ all \ 0<\eta<A.
\end{equation}

b) Assume that $h$ is non-decreasing. Let $A$ be defined by \eqref{def-A}. Then
$$\limsup_{t \to \infty}\frac{1}{t}\log h(t)=A.$$
\end{lemma}

\noindent {\bf Proof:} Part a) is obvious. For part b), note that by Lemma A.1 of \cite{balan-chen16},
$$\bar{A}:=\limsup_{t \to \infty}\frac{1}{t}\log h(t) \leq A.$$
For the other inequality, let $\eta<A$ be arbitrary.
We prove that $\bar{A} \geq \eta$. (The same argument was used for showing that (5.62) implies (5.64) in \cite{FK13}.)
Suppose this is not true. Say $\bar{A}<\eta-\delta$ for some $\delta \in (0,\eta)$. Then $h(t) \leq e^{(\eta-\delta)t}$, for all $t \geq t_0$, for some $t_0$. It follows that
$$\int_0^{\infty}e^{-\eta t}h(t)dt \leq C t_0+\int_{t_0}^{\infty}e^{-\eta t}h(t) dt \leq \int_{t_0}^{\infty}e^{-\eta t} e^{(\eta-\delta)t}dt<\infty,$$
which is a contradiction. $\Box$

\begin{lemma}
\label{lemmaA}
Let $a \in (0,4)$ and $c>0$ be arbitrary. For any $t>0$, define $$A_t=\sum_{n \geq 0}\frac{(ct)^{an}}{\Gamma(an+1)}.$$ Then
$$\lim_{t \to \infty}\frac{1}{t}\log A_t=c.$$
\end{lemma}

\vspace{3mm}

\noindent {\bf Proof:} Note that $A_t=E_a((ct)^a)$, where $E_a(x)=\sum_{n\geq 0} x^n/\Gamma(an+1)$ is the Mittag-Leffler function. If $a \in (0,2)$, we use
Theorem 1.3 (p. 32) of \cite{podlubny99}, and if $a \in [2,4)$, we use Theorem 1.7 (p.35) of \cite{podlubny99}. In both cases, we have:
$$E_a(x)=\frac{1}{a}\exp(x^{1/a})-\frac{x^{-1}}{\Gamma(1-a)}+R(x)$$
where $R(x) \leq Cx^{-2}$ for all $x \geq x_0$, for some $x_0>0$ and $C>0$. Hence
\begin{eqnarray*}
\frac{1}{t}\log A_t &=& -\frac{1}{t} \log a+\frac{1}{t}\log \left(
e^{ct}-\frac{a(ct)^{-1}}{\Gamma(1-a)}+a R(ct) \right)\\
&=& -\frac{1}{t} \log a+c +\frac{1}{t}\log \left(1-\frac{a}{c \Gamma(1-a)}\cdot \frac{1}{t e^{ct}}+ \frac{a R(ct)}{e^{ct}} \right).
\end{eqnarray*}
The conclusion follows. $\Box$

\section{The functionals $\cE$ and $\cE_2$}
\label{appB}

\begin{lemma}
\label{relation-2E-lemma}
Relation \eqref{relation-2E} holds for any non-negative locally integrable function $f$ on $\bR^d$, which is the Fourier transform (in $\cS'(\bR^d)$) of a tempered measure $\mu$ on $\bR^d$ which satisfies \eqref{Dalang-cond}.
\end{lemma}

\noindent {\bf Proof:} Note first that
\begin{equation}
\label{relation-2E-step1}
\cE_2(f) = \sup_{g\in \cF_{2d}}\bigg\{\int_{\bR^d} \int_{\bR^d}
f(\sqrt{2}x_1)g^2(x_1,x_2)dx_1 dx_2
-{1\over 2}\int_{\bR^d} \int_{\bR^d}\vert\nabla g(x_1,x_2)\vert^2dx_1 dx_2\bigg\}=:A.
\end{equation}
This follows considering the one-to-one transformation $g\mapsto \widetilde{g}$ from $\cF_{2d}$ onto itself, where $$\widetilde{g}(y_1,y_2)=g\left(\frac{y_1+y_2}{\sqrt{2}},\frac{y_2-y_1}{\sqrt{2}}\right).$$

We prove that $A=B$, where $B=\cE(f(\sqrt{2} \, \cdot \,))$. First, we show that $B \geq A$. For this, let $g \in \cF_{2d}$ be arbitrary, and define
$\overline{g}(x_1)=\left(\int_{\bR^d}g^2(x_1,x_2)dx_2 \right)^{1/2}$ for all $x_1 \in \bR^d$. Then $\overline{g}\in \cF_d$ and
for any $i=1,\ldots,d$
$$\frac{\partial \overline{g}}{\partial x_1^i}(x_1)=\left(\int_{\bR^d}g^2(x_1,x_2)dx_2 \right)^{-1/2}\int_{\bR^d}g(x_1,x_2)\frac{\partial g}{\partial x_1^i}(x_1,x_2)dx_2.$$
Using the Cauchy-Schwarz inequality and then taking the $dx_1$ integral, we infer that:
$$\int_{\bR^d}\left|\frac{\partial \overline{g}}{\partial x_1^i}(x_1)\right|^2dx_1 \leq \int_{\bR^d} \int_{\bR^d}\left|\frac{\partial g}{\partial x_1^i}(x_1,x_2)\right|^2 dx_2 dx_1.$$ Taking the sum over $i=1,\ldots, d$, we obtain
$\int_{\bR^d} |\nabla \overline{g}(x_1)|^2 dx_1 \leq \int_{\bR^d}\int_{\bR^d}|\nabla g(x_1,x_2)|^2 dx_1 dx_2$.
Hence,
\begin{align*}
B \geq &\int_{\bR^d} f(\sqrt{2}x_1)\overline{g}^2(x_1)dx_1-\frac{1}{2}\int_{\bR^d} |\nabla \bar{g}(x_1)|^2 dx_1\\
\geq &\int_{\bR^d}\int_{\bR^d}f(\sqrt{2}x_1)g^2(x_1,x_2)dx_1 dx_2-\frac{1}{2}\int_{\bR^d}\int_{\bR^d}|\nabla g(x_1,x_2)|^2 dx_1 dx_2.
\end{align*}
The fact that $B \geq A$ follows by taking the supremum over all $g \in \cF_{2d}$.

We now prove that $A \geq B$. Let $g \in \cF_d$ be arbitrary. Define $G(x_1,x_2)=g_(x_1)p_t(x_2)$ where $p_t(x)=(2\pi t)^{-d/2}\exp(-|x|^2/(2t))$ and $t>0$. Then $G \in \cF_{2d}$. Denote $x_1=(x_1^1,\ldots,x_1^d)$ and $x_2=(x_2^1,\ldots,x_2^d)$. For any $i=1,\ldots,d$,
$$\frac{\partial G}{\partial x_1^i}(x_1,x_2)=\frac{\partial g}{\partial x_1^i}(x_1)p_t(x_2) \quad \mbox{and} \quad \frac{\partial G}{\partial x_2^i}(x_1,x_2)=-\frac{1}{2t}g(x_1)p_t(x_2)^{1/2}x_2^i.$$
From this, we obtain:
$$\int_{\bR^d}\int_{\bR^d}|\nabla G(x_1,x_2)|^2 dx_1 dx_2=\int_{\bR^d}|\nabla g(x_1)|^2 dx_1+\frac{1}{4t}c,$$
where $c=E|Z|^2$ and $Z$ is a $N_d(0,I)$-random vector. Hence
\begin{align*}
 A\geq &\int_{\bR^d} \int_{\bR^d} f(\sqrt{2}x_1)G^2(x_1,x_2)dx_1dx_2-\frac{1}{2}\int_{\bR^d}\int_{\bR^d} |\nabla G(x_1,x_2)|^2 dx_1 dx_2\\
 =&\int_{\bR^d}f(\sqrt{2}x_1)g^2(x_1)dx_1 -\frac{1}{2}\int_{\bR^d}|\nabla g(x_1)|^2 dx_1 -\frac{1}{8t}c.
\end{align*}
We let $t \to \infty$. Then we take the supremum over all $g \in \cF_d$. $\Box$

\vspace{3mm}

The next result gives a scaling property of the functionals $\cE(f)$ and $\cE_2(f)$.

\begin{lemma}
\label{scaling-E}
Let $f$ be a non-negative locally integrable function on $\bR^d$ which satisfies \eqref{scaling-f} for some $\alpha>0$. Then for any $\theta>0$,
$$\cE(\theta f)= \theta^{\frac{2}{2-\alpha}}\cE(f) \quad \mbox{and} \quad \cE_2(\theta f)= \theta^{\frac{2}{2-\alpha}}\cE_2(f).$$
Moreover, $\cE_2(f)=2^{-\frac{\alpha}{2-\alpha}}\cE(f)$.
\end{lemma}

\noindent {\bf Proof:} Due to \eqref{relation-2E}, we only need to prove the scaling property of $\cE(f)$. This can be verified using the one-to-one transformation $g \mapsto \widetilde{g}$ from $\cF_d$ onto itself, where
$\tilde g(y)=\theta^{-\frac{d}{4-2\alpha}}g(\theta^{-\frac1{2-\alpha}}y)$. Indeed, noting that $|\nabla \widetilde g(y)|^2 = \theta^{-\frac{d+2}{2-\alpha}} |\nabla g(\theta^{-\frac1{2-\alpha}}y)|^2$, we have
$$\theta \int_{\bR^d}f(x)g^2(x)dx=\theta^{\frac{2}{2-\alpha}} \int_{\bR^d}f(y)\widetilde{g}^2(y)dy, \quad \int_{\bR^d}|\nabla g(x)|^2dx=\theta^{\frac{2}{2-\alpha}}\int_{\bR^d}|\nabla \widetilde{g}(y)|^2dy,$$
and
\begin{align*}
\cE(\theta f)=&\sup_{g\in {\cal F}_d}\bigg\{\theta\int_{\bR^d}
f(x)g^2(x)dx
-{1\over 2}\int_{\bR^d}\vert\nabla g(x)\vert^2dx\bigg\}\\
=&\sup_{g\in {\cal F}_d}\bigg\{\theta^{\frac{2}{2-\alpha}}\int_{\bR^d}
f(y)\widetilde{g}^2(y)dy
-\theta^{\frac{2}{2-\alpha}}{1\over 2}\int_{\bR^d}|\nabla \widetilde{g}(y)|^2dy\bigg\}=\theta^{\frac{2}{2-\alpha}} \cE(f).
\end{align*}

For the last equality, note that by \eqref{relation-2E}, \eqref{scaling-f} and the scaling property of $\cE$, we have:
$$\cE_2(f)=\cE(f(\sqrt{2} \, \cdot \,))=\cE(2^{-\alpha/2}f)=2^{-\frac{\alpha}{2-\alpha}}\cE(f).$$
$\Box$

\vspace{3mm}
 {
To prove the statement of Remark \ref{comparison-remark}, we need the following result.}

\begin{lemma}
\label{E-A-functional}
 {If $f$ satisfies \eqref{scaling-f}, then for any $A>0$,
$$\cE(f,A)=A^{\frac{\alpha}{\alpha-2}}\cE(f,1),$$
where
\begin{equation}
\label{def-E-A}
\cE(f,A)= \sup_{g\in \cF_d}\bigg\{\int_{\bR^d}
f(x)g^2(x)dx
-A\int_{\bR^d}\vert\nabla g(x)\vert^2dx\bigg\}.
\end{equation}
}
\end{lemma}

\noindent  {{\bf Proof:} As in the proof of Lemma \ref{scaling-E}, we consider the one-to-one transformation of $\cF_d$ onto itself, where $\tilde g(y)=A^{-\frac{d}{4-2\alpha}}g(A^{-\frac1{2-\alpha}}y)$. Then
$$\int_{\bR^d}f(y)\widetilde{g}^2(y)dy=A^{\frac{\alpha}{\alpha-2}}\int_{\bR^d}f(x)g^2(x)dx , \quad A\int_{\bR^d}|\nabla \widetilde{g}(y)|^2dy=A^{\frac{\alpha}{\alpha-2}}\int_{\bR^d}|\nabla g(x)|^2dx,$$
and
\begin{align*}
A^{\frac{\alpha}{\alpha-2}}\cE(f,1)=&\sup_{g\in {\cal F}_d}\bigg\{A^{\frac{\alpha}{\alpha-2}}\int_{\bR^d}
f(x)g^2(x)dx
-A^{\frac{\alpha}{\alpha-2}}\int_{\bR^d}\vert\nabla g(x)\vert^2dx\bigg\}\\
=&\sup_{g\in {\cal F}_d}\bigg\{\int_{\bR^d}
f(y)\widetilde{g}^2(y)dy
-A\int_{\bR^d}|\nabla \widetilde{g}(y)|^2dy\bigg\}=\cE(f,A).
\end{align*}
$\Box$
}

When $d=1$, the function $f(x)=|x|^{2H-2}$ with $1/2<H<1$ satisfies the scaling property \eqref{scaling-f} with $\alpha=2-2H$, and hence, $\cE(\theta f)=\theta^{1/H}\cE(f)$ by Lemma \ref{scaling-E}. The following result shows that this property continues to hold in the case $H<1/2$.

\begin{lemma}
\label{scaling-E-Gamma}
If $\Gamma$ is the distribution given by \eqref{def-Gamma}, then for any $\theta>0$,
$$\cE(\theta \Gamma)= \theta^{1/H}\cE(\Gamma) \quad \mbox{and} \quad \cE_2(\theta \Gamma)= \theta^{1/H}\cE_2(\Gamma).$$
Moreover, $\cE_2(\Gamma)=2^{-\frac{1-H}{H}}\cE(\Gamma)$.
\end{lemma}

\noindent {\bf Proof:} We use the same arguments as in the proof of Lemma \ref{scaling-E} (with $d=1$ and $\alpha=2-2H$), but we express all the quantities in the Fourier mode. To show the scaling property of $\cE(\Gamma)$, it suffices to consider the one-to-one transformation $h \mapsto \widetilde{h}$ from $\cA_1$ onto itself, where $\widetilde{h}(\eta)=\theta^{\frac{1}{4H}}h(\theta^{\frac{1}{2H}}\eta)$. Then
$$\theta \int_{\bR} (h*h)(\xi)\mu(d\xi)=\theta^{\frac{1}{H}}\int_{\bR}(\widetilde{h}*\widetilde{h})
(\eta)\mu(d\eta), \quad \int_{\bR}|\xi|^2 |h(\xi)|^2 d\xi=\theta^{\frac{1}{H}}\int_{\bR}
|\eta|^2 |\widetilde{h}(\eta)|^2d\eta.$$
The scaling property of $\cE_2(\Gamma)$ will follow from the last relation.

To prove the last relation, we note that
\begin{align*}
\cE_2(\Gamma)&=\sup_{h \in \cA_2}
\bigg\{\int_{\bR}
 (h*h)(\sqrt{2}\xi,0)\mu(d\xi)-\frac{1}{2}\int_{\bR}\int_{\bR}(|\xi_1|^2+|\xi_2|^2) |h(\xi_1,\xi_2)|^2 d\xi_1 d\xi_2  \bigg\}\\
 &=\sup_{h \in \cA_1} \bigg\{\int_{\bR}
 (h*h)(\sqrt{2}\xi)\mu(d\xi)-\frac{1}{2}\int_{\bR}|\xi|^2 |h(\xi)|^2 d\xi  \bigg\}\\
 &=\cE(2^{H-1}\Gamma)=2^{\frac{H-1}{H}}\cE(\Gamma).
\end{align*}
The first equality above is the analogue of \eqref{relation-2E-step1} in Fourier mode, and can be proved using the one-to-one transformation $h \mapsto \widetilde{h}$ from $\cA_2$ onto itself, where $\widetilde{h}(\eta_1,\eta_2)=h(\frac{\eta_1+\eta_2}{\sqrt{2}}, \frac{\eta_2-\eta_1}{\sqrt{2}})$. The second equality is the analogue of the equality  $A=B$ in the proof of Lemma \ref{relation-2E-lemma}, and can be proved by similar methods. The third equality follows by a change of variables, and the last equality is due to the scaling property of $\cE(\Gamma)$. $\Box$

\end{document}